\documentclass[a4paper,12pt]{article}
\usepackage{amsmath, amsfonts, amssymb}

\makeatletter
\newbox\bk@bxb
\newbox\bk@bxa
\newif\if@bkcont
\newcount\bk@lcnt

\def\breakboxskip{2pt}
\def\breakboxparindent{1.8em}

\def\breakbox{\vskip\breakboxskip\relax
\setbox\bk@bxb\vbox\bgroup
\advance\linewidth -2\fboxrule
\hsize\linewidth\@parboxrestore
\parindent\breakboxparindent\relax}

\def\bk@split{%
\@tempdimb\ht\bk@bxb 
\advance\@tempdimb\dp\bk@bxb
\setbox\bk@bxa\vsplit\bk@bxb to\z@ 
\setbox\bk@bxa\vbox{\unvbox\bk@bxa}
\setbox\@tempboxa\vbox{\copy\bk@bxa\copy\bk@bxb}
\advance\@tempdimb-\ht\@tempboxa
\advance\@tempdimb-\dp\@tempboxa}

\def\bk@addfsepht{%
\setbox\bk@bxa\vbox{\vskip\fboxsep\box\bk@bxa}}

\def\bk@addskipht{%
\setbox\bk@bxa\vbox{\vskip\@tempdimb\box\bk@bxa}}

\def\bk@addfsepdp{%
\@tempdima\dp\bk@bxa
\advance\@tempdima\fboxsep
\dp\bk@bxa\@tempdima}

\def\bk@addskipdp{%
\@tempdima\dp\bk@bxa
\advance\@tempdima\@tempdimb
\dp\bk@bxa\@tempdima}

\def\bk@line{%
\hbox to \linewidth{%
\hskip-2\fboxsep\vrule \@width\fboxrule\hskip.5\fboxsep\vrule \@width\fboxrule\hskip1.5\fboxsep
\box\bk@bxa\hfil
}}%

\def\endbreakbox{\egroup
\ifhmode\par\fi{\noindent\bk@lcnt\@ne
\@bkconttrue\baselineskip\z@\lineskiplimit\z@
\lineskip\z@\vfuzz\maxdimen
\bk@split\bk@addfsepht\bk@addskipdp
\ifvoid\bk@bxb 
\def\bk@fstln{\bk@addfsepdp
\hskip-\parindent\vbox{\llap{\raisebox{-2ex}{\rule{1.5\fboxsep}{\fboxrule}\hskip.5\fboxsep}}\bk@line\llap{\rule{1.5\fboxsep}{\fboxrule}\hskip.5\fboxsep}}}

\else 
\def\bk@fstln{\vbox{\llap{\raisebox{-2ex}{\rule{1.5\fboxsep}{\fboxrule}\hskip.5\fboxsep}}\bk@line}\hfil%
\advance\bk@lcnt\@ne
\loop
\bk@split\bk@addskipdp\leavevmode
\ifvoid\bk@bxb 
\@bkcontfalse\bk@addfsepdp
\vtop{\bk@line\llap{\rule{2\fboxsep}{\fboxrule}}}%

\else 
\bk@line
\fi
\hfil\advance\bk@lcnt\@ne
\if@bkcont\repeat}%
\fi
\leavevmode\bk@fstln\par}\vskip\breakboxskip\relax}

\newcommand{\fracb}[2]{\frac{\raisebox{-.7ex}{$\displaystyle #1$}}{\raisebox{.7ex}{$\displaystyle #2$}}}

\newcommand{\scalv}[2]{\langle#1,#2\rangle}
\newcommand{\scal}[2]{\langle#1\mid#2\rangle}

\def\smp{\smallskip\par}
\def\un{{\bf 1}}
\def\zero{\{0\}}
\def\pf{\noindent{\bf Proof~:}\ }

\def\findemo{~\leaders\hbox to 1em{\hss\  \hss}\hfill~\raisebox{.5ex}{\framebox[1ex]{}}\smp}
\def\spn{\bigskip\par\noindent}
\def\mpn{\medskip\par\noindent}
\def\smpn{\smallskip\par\noindent}
\def\normal{\mathop{\trianglelefteq}}
\def\sp{\bigskip\par}
\def\smp{\smallskip\par}
\def\smpn{\smallskip\par\noindent}

\def\mpoint{\;\;.}
\def\mvirg{\;\;,}

\def\Res{{\rm Res}}
\def\Ind{{\rm Ind}}

\def\Hom{{\rm Hom}}
\def\End{{\rm End}}
\def\Ext{{\rm Ext}}

\def\Out{{\rm Out}}
\def\Aut{{\rm Aut}}

\def\Ker{{\rm Ker}}

\def\Id{{\rm Id}}

\def\op{^{op}}

\def\Z{\mathbb{Z}}

\def\F{\mathbb{F}}
\def\Q{\mathbb{Q}}

\def\C{\mathbb{C}}

\newcommand{\dirsum}[1]{\mathop{\oplus}_{#1}\limits}
\newcommand{\romain}[1]{\uppercase\expandafter{\romannumeral #1}}

\newcommand{\flh}[2]{\mathop{\hbox to 12mm{\rightarrowfill}}_{\displaystyle #2}^{\displaystyle #1}\limits}
\newcommand{\sflh}[2]{\mathop{\hbox to 12mm{\rightarrowfill}}_{\scriptstyle #2}^{\scriptstyle #1}\limits}

\newcommand{\gMod}[1]{#1{\hbox{-}\mathsf{Mod}}}

\newcommand{\sur}[1]{\,\overline{\! #1}}
\newcommand{\sumb}[2]{\mathop{\sum}_{{\scriptstyle #1}\atop {\scriptstyle #2}}}

\def\op{^{op}}

\newcommand{\carre}[8]{\begin{array}{ccc}
#1&\mathop{\hbox to 12mm{\rightarrowfill}}^{\displaystyle{#2}}\limits&#3\\
\llap{$\displaystyle{#4}$}\left\downarrow\vbox to 6mm{}\right. & & \left\downarrow\vbox to 6mm{}\right.\rlap{$\displaystyle{#5}$}\\
#6&\mathop{\hbox to 12mm{\rightarrowfill}}_{\displaystyle #7}\limits&#8\\
\end{array}}

\newcommand{\carrem}[8]{\begin{array}{ccc}
#1&\mathop{\hbox to 12mm{\rightarrowfill}}^{\displaystyle #2}\limits&#3\\
\llap{$\displaystyle #4$}\left\uparrow\vbox to 6mm{}\right. & & \left\uparrow\vbox to 6mm{}\right.\rlap{$\displaystyle #5$}\\
#6&\mathop{\hbox to 12mm{\rightarrowfill}}_{\displaystyle #7}\limits&#8\\
\end{array}}

\newenvironment{enonce}[1]{\pagebreak[2]\refstepcounter{subsection}\refstepcounter{prop}\smpn{{\bf \thesection.\arabic{prop}.\ #1:}}\begin{it} }{\end{it}\smp}
\newenvironment{enonce*}[1]{\pagebreak[2]\smpn{#1:}\begin{it} }{\end{it}\smp}
\newcommand{\result}[1]{\begin{enonce}{#1}}
\def\fresult{\end{enonce}}
\newcommand{\npar}{\smallskip\par\noindent\pagebreak[2]\refstepcounter{subsection}\refstepcounter{prop}{\bf \thesection.\arabic{prop}.\ \ }}

\newenvironment{mth}[1]{\begin{breakbox}\begin{enonce}{#1}}{\end{enonce}\end{breakbox}}
\newenvironment{mth*}[1]{\begin{breakbox}\begin{enonce*}{#1}}{\end{enonce*}\end{breakbox}}
\newenvironment{rem}[1]{\refstepcounter{subsection}\refstepcounter{prop} \mpn{{\bf \thesection.\arabic{prop}.}\ \bf#1:}}{\smp}

\def\dom{\backslash}
\makeatletter
\renewenvironment{enumerate}{\ifnum \@enumdepth >3 \@toodeep\else
      \advance\@enumdepth \@ne
      \edef\@enumctr{enum\romannumeral\the\@enumdepth}\list
      {\csname label\@enumctr\endcsname}{\setlength{\topsep}{1ex}\setlength{\itemsep}{0pt}\usecounter
        {\@enumctr}\def\makelabel##1{\hss\llap{##1}}}\fi}{\endlist}
\renewenvironment{itemize}{\ifnum \@itemdepth >3 \@toodeep\else \advance\@itemdepth \@ne
\edef\@itemitem{labelitem\romannumeral\the\@itemdepth}%
\list{\csname\@itemitem\endcsname}{\setlength{\topsep}{1ex}\setlength{\itemsep}{0pt}\def\makelabel##1{\hss\llap{##1}}}\fi}
{\endlist}
\makeatother
\makeatletter
\def\@sect#1#2#3#4#5#6[#7]#8{\ifnum #2>\c@secnumdepth
    \let\@svsec\@empty\else
    \refstepcounter{#1}\edef\@svsec{\csname the#1\endcsname .\hskip .5em}\fi
    \@tempskipa #5\relax
     \ifdim \@tempskipa>\z@
       \begingroup #6\relax
         \@hangfrom{\hskip #3\relax\@svsec}{\interlinepenalty \@M #8\par}%
       \endgroup
      \csname #1mark\endcsname{#7}\addcontentsline
        {toc}{#1}{\ifnum #2>\c@secnumdepth \else
                     \protect\numberline{\csname the#1\endcsname}\fi
                   #7}\else
       \def\@svsechd{#6\hskip #3\relax  
                  \@svsec #8\csname #1mark\endcsname
                     {#7}\addcontentsline
                          {toc}{#1}{\ifnum #2>\c@secnumdepth \else
                            \protect\numberline{\csname the#1\endcsname}\fi
                      #7}}\fi
    \@xsect{#5}}
\def\section{\@startsection {section}{1}{\z@}{-3.5ex plus-1ex minus
    -.2ex}{2.3ex plus.2ex}{\center\reset@font\large\bf}}  

\makeatother
\renewenvironment{equation}{\refstepcounter{subsection}\refstepcounter{prop}$$}{\leqno{\bf (\theprop)}$$}

\def\mar[#1]{\ar@{-}[#1]|-{\object@{<}}}
\def\marb[#1]{\ar@{-}[#1]|{\object+{  }}}

\def\mpoint{\;.}
\def\mvirg{\;,}
\def\CC{\mathcal{C}}
\newcommand{\ncalv}[2]{\scalv{#1}{#2}^{\natural}}
\renewcommand{\fracb}[2]{\frac{\scriptstyle#1}{\scriptstyle#2}}
\begin{document}
\centerline{\Large\bf Some simple biset functors}\vspace{.5cm}\par
\centerline{\bf Serge Bouc}\vspace{1cm}\par
{\footnotesize {\bf Abstract:} Let $p$ be a prime number, let $H$ be a finite $p$-group, and let $\F$ be a field of characteristic 0, considered as a trivial $\F \Out(H)$-module. The main result of this paper gives the dimension of the evaluation $S_{H,\F}(G)$ of the simple biset functor $S_{H,\F}$ at an arbitrary finite group~$G$. A closely related result is proved in the last section: for each prime number $p$, a Green biset functor $E_p$ is introduced, as a specific quotient of the Burnside functor, and it is shown that the evaluation $E_p(G)$ is a free abelian group of rank equal to the number of conjugacy classes of $p$-elementary subgroups of $G$.
}\vspace{2ex}\par
\vspace{2ex}\par
\vspace{1ex}
\section{Introduction}
Let $R$ be a commutative ring. The {\em biset category} $R\CC$ over $R$ has finite groups as objects, with morphisms $\Hom_{R\CC}(G,H)=R\otimes_\Z B(H,G)$, where $B(H,G)$ is the Burnside group of $(H,G)$-bisets. The composition of morphisms is induced by the usual tensor product of bisets. A {\em biset functor} over $R$ is an $R$-linear functor from $R\CC$ to the category $\gMod{R}$ of $R$-modules. Biset functors over $R$ form an abelian category, where morphisms are natural transformations of functors. They have proved a useful tool in various aspects of the representation theory of finite groups (see \cite{both}, \cite{dadegroup}, \cite{burnsideunits}, \cite{bisetfunctorsMSC}), and they are still the object of active research (\cite{garcia}, \cite{barker-ogut}, \cite{centros}, \cite{rognerud-evaluations}, \cite{coskun-yalcin}, \cite{coskun-yilmaz}, \cite{boltje-raggi-valero}, \cite{boltje-coskun}, \cite{barsotti}, \ldots). \par
The simple biset functors over $R$ are parametrized (\cite{doublact}, Proposition 2) by equivalence classes of pairs $(H,W)$, where $H$ is a finite group, and $W$ is a simple $R\Out(H)$-module - the simple functor parametrized by $(H,W)$ being denoted $S_{H,W}$. However for a finite group $G$, the computation of the evaluation $S_{H,W}(G)$ is generally quite hard: in Theorem~4.3.20 of~\cite{bisetfunctorsMSC}, this evaluation is shown to be equal to the image of a complicated linear map. Assuming that $R$ is a field - which is always possible when dealing with simple functors - the dimension of $S_{H,W}(G)$ is given by Theorem~7.1 of~\cite{boustathe}, as the rank of a yet complicated bilinear form with values in $R$. \par
Let $\F$ be a field of characteristic 0, let $p$ be a prime number, and $H$ be a finite $p$-group. The present paper is mainly devoted to the computation of the dimension of the evaluation $S_{H,\F}(G)$, where $G$ is an arbitrary finite group, and $\F$ is the trivial $\F\Out(H)$-module. The result is as follows:
\pagebreak[3]
\begin{mth*}{{\bf Theorem}} Let $\F$ be a field of characteristic 0, let $p$ be a prime number, and $H$ be a finite $p$-group. Let moreover $G$ be a finite group.
\begin{enumerate}
\item If $H=\un$, the dimension of $S_{H,\F}(G)$ is equal to the number of conjugacy classes of cyclic subgroups of $G$.
\item If $H\cong C_p\times C_p$, the dimension of $S_{H,\F}(G)$ is equal to the number of conjugacy classes of non-cyclic $p$-elementary subgroups of $G$.
\item If $H$ is any other finite $p$-group, the dimension of $S_{H,\F}(G)$ is equal to the number of conjugacy classes of sections $(T,S)$ of $G$ such that $T/S\cong H$ and $T$ is $p$-elementary.
\end{enumerate}
\end{mth*}
In the last section of this paper, for each prime number $p$, we introduce a Green biset functor $E_p$, closely related to the two first assertions of the above theorem. Green biset functors have been defined in~\cite{bisetfunctorsMSC}, Section~8.5. They are ring objects in the category of biset functors. For a finite group $G$, we denote by $F_p(G)$ the set of elements of the Burnside group $B(G)$ which vanish when restricted to all $p$-elementary subgroups of $G$, and we show that this actually defines a biset subfunctor $F_p$ of $B$. The functor $E_p$ is defined as the quotient $B/F_p$, and it then inherits from $B$ a Green biset functor structure (over $\Z$). We show moreover that its evaluation $E_p(G)$ at a finite group $G$ is a free abelian group of rank equal to the number of conjugacy classes of $p$-elementary subgroups of $G$. We also show that the biset functor $\F E_p=\F\otimes_\Z E_p$ fits in a non split short exact sequence
$$0\to S_{(C_p)^2,\F}\to \F E_p\to S_{\un,\F}\to 0$$
of biset functors over $\F$. In the case $\F=\Q$, the restriction of this sequence to $p$-groups is the short exact sequence of Theorem D of~\cite{both}, involving the Dade functor $\Q D$, the Burnside functor $\Q B$, and the functor of rational representation $\Q R_\Q$.

\section{Preliminary results}
Recall (Sections~5.3 and~5.4 of~\cite{bisetfunctorsMSC}) that for a normal subgroup $N$ of a finite group $G$, the rational number $m_{G,N}$ is defined by
$$m_{G,N}=\fracb{1}{|G|}\sum_{\substack{X\leq G\\\rule{0ex}{1.5ex}XN=G}}|X|\mu(X,G)\mvirg$$
where $\mu$ is the M\"obius function of the poset of subgroups of $G$. The group~$G$ is called a {\em $B$-group} if $m_{G,N}=0$ for any non-trivial normal subgroup $N$ of~$G$. Any finite group $G$ has a largest quotient $B$-group $\beta(G)$, unique up to isomorphism. If $N\normal G$, then $m_{G,N}=0$ if and only if $\beta(G)\cong\beta(G/N)$.
\begin{mth}{Lemma}{\rm [M. Baumann~\cite{baumann-jofa} - See also \cite{doublact}, 8), p 713] } \label{EL}Let $L$ be a finite group, let $p$ be a prime, and let $E$ be an elementary abelian $p$-group on which $L$ acts irreducibly, faithfully, and such that $H^1(L,E)=\zero$. Then the group $G=E\rtimes L$ is a $B$-group.
\end{mth}
\pf First as $E$ is $L$-simple, it follows that $E$ is a minimal normal subgroup of $G$. Let $N$ be any normal subgroup of $G$. Then $N\cap E$ is equal to $E$ or $\un$. So if $N\ngeq E$, then $N\cap E=\un$, and $N$ centralizes $E$. But $C_G(E)=E\cdot C_L(E)=E$, since $L$ acts faithfully on $E$. Thus $N\leq E$, hence $N=\un$.\par
It follows that $E$ is the unique minimal normal subgroup of $G$. By Proposition~5.6.4 of~\cite{bisetfunctorsMSC}, since $E$ is abelian, 
\begin{equation}\label{minimal normal}
m_{G,E}=1-\frac{|K_G(E)|}{|E|}\mvirg
\end{equation}
where $K_G(E)$ is the set of complements of $E$ in $G$. The group $E$ acts by conjugation on $K_G(E)$, and the normalizer in $E$ of $K\in K_G(N)$ is equal to the group $E^K$ of fixed points of $K$ on $E$. Since $E$ is $K$-simple, and $K$-faithful, this is equal to $\un$. Thus $E$ acts freely on $K_G(E)$. Since $H^1(L,E)=\zero$, the set $K_G(E)$ is a single conjugacy class, i.e. a single $E$-orbit. Thus $|K_G(E)|=|E|$, and $m_{G,E}=0$. It follows that $G$ is a $B$-group.\findemo
\npar Recall that a finite group $G$ is called {\em cyclic modulo a prime number $p$} if $G/O_p(G)$ is cyclic, and that $G$ is called {\em $p$-elementary} if $G\cong P\times C$, where $P$ is a $p$-group and $C$ is a cyclic group.
\begin{mth}{Lemma} \label{p-elementary}Let $p$ be a prime number, and $G$ be a finite group.
\begin{enumerate}
\item {\rm [M. Baumann~\cite{baumann-jofa}]} The group $\beta(G)$ is cyclic modulo $p$ if and only if $G$ is cyclic modulo $p$.
\item The group $\beta(G)$ is a $p$-group if and only if $G$ is $p$-elementary.
\end{enumerate}
\end{mth}
\pf For Assertion~1, use the fact that by a theorem of Conlon, the subspace $NC_p(G)$ of $\Q B(G)$ generated by the idempotents $e_H^G$, where $H$ is {\em not} cyclic modulo $p$, is equal to the kernel of the morphism $\Q B(G)\to \Q pp_k(G)$. In particular, the correspondence $G\mapsto NC_p(G)$ is a biset subfunctor of $\Q B$. It follows that there exists a family $\mathcal{B}$ of $B$-groups such that for any group~$G$, the space $NC_p(G)$ is the $\Q$-vector subspace of $\Q B(G)$ generated by the idempotents $e_H^G$, where $\beta(H)\in\mathcal{B}$. The family $\mathcal{B}$ consists of those $B$-groups $H$ for which $e_H^H\in NC_p(H)$, i.e. the $B$-groups which are not cyclic modulo $p$. Now for any group $G$, and any subgroup $H$ of $G$, the idempotent $e_H^G$ is in $NC_p(G)$ if and only if $\beta(G)\in\mathcal{B}$, on the one hand, but also if and only if $H$ is not cyclic modulo $p$. Hence $\beta(H)$ is not cyclic modulo $p$ if and only if $H$ is not cyclic modulo $p$. This proves Assertion~1.\par
For Assertion~2, clearly, if $G$-is $p$-elementary, one can assume $G\cong P\times C$, where $P$ is a $p$-group, and $C$ is a cyclic $p'$-group. By Proposition~5.6.6 of~\cite{bisetfunctorsMSC}, this implies $\beta(G)\cong\beta(P)\times \beta(C)\cong \beta(P)$, since $\beta(C)=\un$. Hence $\beta(G)$ is a $p$-group.\par
Conversely, suppose that $\beta(G)$ is a $p$-group. In particular, it is cyclic modulo $p$, hence $G$ is cyclic modulo $p$, by Lemma~\ref{EL}.
The Frattini subgroup $\Phi(P)$ of $P$ is a normal subgroup of $G$, and $G/\Phi(P)\cong \sur{P}\rtimes C$, where $\sur{P}$ is the elementary abelian group $P/\Phi(P)$. Suppose that the $\F_pC$-module $\sur{P}$ admits a simple quotient $E$ with non-trivial $C$-action (that is, not isomorphic to $\F_p$). Then the action of $C$ on $E$ has a kernel $D<C$, and the group $E\rtimes (C/D)$ is a quotient of $\sur{P}\rtimes C$, hence a quotient of $G$. But $E\rtimes (C/D)$ is a $B$-group by Lemma~\ref{EL}: Indeed $E$ is $(C/D)$-simple and faithful by construction, and $H^1(C/D,E)=\zero$, since $C/D$ is a $p'$-group. \par
Now $E\rtimes (C/D)$ is a $B$-group, which is not a $p$-group, since $D\neq C$, and it is a quotient of $G$, hence of $\beta(G)$, which is a $p$-group. This is a contradiction.\par
Hence $C$ acts trivially on $\sur{P}$. But for any $p$-group $P$, the kernel of the morphism $\Aut(P)\to\Aut\big(P/\Phi(P)\big)$ is a $p$-group. As $C$ is a $p'$-group, and acts trivially on $\sur{P}$, it acts trivially on $P$. Thus $G\cong P\times C$, as was to be shown.\findemo
\begin{mth}{Lemma} \label{complements}Let $p$ be a prime number, and $P$ be a finite $p$-group.
\begin{enumerate}
\item Let $Q$ be a normal subgroup of $P$. Then $Q\cap \Phi(P)=\un$ if and only if $Q$ is elementary abelian and central in $P$, and admits a complement in~$P$.
\item Let $Q$ and $R$ be normal subgroups of $P$, such that $|Q|=|R|$. Then $Q\cap \Phi(P)=\un=R\cap\Phi(P)$ if and only if $Q$ and $R$ are elementary abelian and central in $P$, and admit a common complement in $P$. \par
In this case, set $H=P/R\cong P/Q$, and denote by $\gamma$ the rank of the group $H/\Phi(H)$. If $Q$ and $R$ have rank $m$, and if $Q\cap R\Phi(P)$ has rank $m-s$,  the number of common complements of $Q$ and $R$ in $P$ is equal to
$$(p^s-1)(p^{s-1}-1)\cdots(p-1)p^{\binom{s}{2}+s(m-s)+m(\gamma-s)}\mpoint$$
\end{enumerate}
\end{mth}
\pf For Assertion~1, if $Q$ is elementary abelian and central in $P$, and admits a complement $L$, then $P=Q\times L$. Thus $\Phi(P)=\un\times \Phi(L)$, hence $Q\cap\Phi(P)=\un$. Conversely, if $Q\cap\Phi(P)=\un$, then $Q$ maps injectively into $P/\Phi(P)$, so $Q$ is elementary abelian. Let $L\geq \Phi(P)$ be a subgroup of~$P$ such that $L/\Phi(P)$ is a complement of $Q\Phi(P)/\Phi(P)$ in the $\F_p$-vector space $P/\Phi(P)$. Then $Q\Phi(P)L=P$, thus $QL=P$, and $Q\Phi(P)\cap L=\Phi(P)$, i.e. $Q\cap L\leq Q\cap\Phi(P)=\un$. Since $L\geq \Phi(P)$, it follows that $L\normal P$, thus $[L,Q]\leq L\cap Q=\un$, and $Q$ is central in $P$. \par
For Assertion 2, let $Q$ and $R$ be normal subgroups of~$P$ with $|Q|=|R|$. If $Q$ and $R$ are elementary abelian central subgroups of $P$ with a common complement in $P$, then $Q\cap\Phi(P)=R\cap\Phi(P)=\un$ by Assertion~1. Conversely, if $Q\cap \Phi(P)=\un$ and $R\cap \Phi(P)=\un$, then $Q$ and $R$ are elementary abelian and central in $P$ by Assertion~1. If $L$ is a complement of $Q$ in $P$, then $P=Q\times L$ for $Q$ is central in $P$, thus $L\normal P$, and $P/L\cong Q$ is elementary abelian. Thus $L\geq\Phi(P)$, and $L/\Phi(P)$ is a complement of $Q\Phi(P)/\Phi(P)$ in $P/\Phi(P)$. Conversely if $L/\Phi(P)$ is a complement of $Q\Phi(P)/\Phi(P)$ in $P/\Phi(P)$, then $L$ is a complement of $Q$ in $P$, by the argument used in the proof of Assertion~1.\par
So finding a common complement to $Q$ and $R$ in $P$ amounts to finding a common complement of $\widetilde{Q}=Q\Phi(P)/\Phi(P)$ and $\widetilde{R}=R\Phi(P)/\Phi(P)$ in $\widetilde{P}=P/\Phi(P)$. Moreover $|\widetilde{Q}|=|Q|=|R|=|\widetilde{R}|$. The $\F_p$-vector space $\widetilde{P}$ can be split as $\widetilde{P}=I\oplus E\oplus F\oplus V$, where $I=\widetilde{Q}\cap\widetilde{R}$, where $E$ is a complement of $I$ in $\widetilde{Q}$ and $F$ is a complement of $I$ in $\widetilde{R}$, and $V$ is a complement of $\widetilde{Q}+\widetilde{R}$ in $\widetilde{P}$. Then $L=F\oplus V$ is a complement of $\widetilde{Q}$ in $\widetilde{P}$, and all the other complements of $\widetilde{Q}$ are of the form $\{\big(\varphi(x),x\big)\mid x\in L\}$, where $\varphi:L\to \widetilde{Q}$ is a group homomorphism. In other words, any complement $L'$ of $\widetilde{Q}$ is of the form
$$L'=\{\big(a(f)+b(v),c(f)+d(v),f,v\big)\mid f\in F,\,v\in V\}\mvirg$$
where $a:F\to I$, $b:V\to I$, $c:F\to E$ and $d:V\to E$ are group homomorphisms. The group $L'$ is a complement of $\widetilde{R}$ if and only if its intersection with $\widetilde{R}$ is trivial, or equivalently if $c$ is injective, hence an isomorphism, since $|E|=|F|$.\par
It follows that the number of common complements of $\widetilde{Q}$ and $\widetilde{R}$ in $\widetilde{P}$ is equal to the number of 4-tuples $(a,b,c,d)$, where $c$ is an isomorphism. Hence
\begin{eqnarray*}
|K_P(Q)\cap K_P(R)|&=&|\Aut(E)||\Hom(F,I)||\Hom(V,I)||\Hom(V,E)|\\
&=&|\Aut(E)||\Hom(F,I)||\Hom(V,\widetilde{Q})|\mpoint
\end{eqnarray*}
Moreover
$$I=\big(Q\Phi(P)\cap R\Phi(P)\big)/\Phi(P)=\big(Q\cap R\Phi(P)\big)\Phi(P)/\Phi(P)\cong Q\cap R\Phi(P)$$
has rank $m-s$, and $F\cong E\cong \widetilde{Q}/I$ has rank $s$.  Finally 
$$V\cong \widetilde{P}/(\widetilde{Q}\widetilde{R})\cong \big(P/R\Phi(P)\big)\big/\big(QR\Phi(P)/R\Phi(P)\big)$$ 
has rank $\gamma-s$, since $P/R\Phi(P)\cong H/\Phi(H)$, as $\Phi(P/R)=R\Phi(P)/R$, and since $QR\Phi(P)/R\Phi(P)\cong Q/\big(Q\cap R\Phi(P)\big)$. This completes the proof.\findemo
\begin{mth}{Corollary} \label{ExH} Let $P$ be a finite $p$-group, and $M$ be a normal subgroup of~$P$. Then
$$P/\big(M\cap \Phi(P)\big)\cong E\times (P/M)\mvirg$$
where $E=M/\big(M\cap \Phi(P)\big)$ is elementary abelian.
\end{mth}
\pf The normal subgroup $\sur{M}=M/\big(M\cap \Phi(P)\big)$ of $\sur{P}=P/\big(M\cap \Phi(P)\big)$ intersect the Frattini subgroup $\Phi(\sur{P})=\Phi(P)/\big(M\cap \Phi(P)\big)$ trivially, hence there exists a subgroup $L$ of $\sur{P}$ such that $\sur{P}=\sur{M}\times L$. Moreover $L\cong\sur{P}/\sur{M}\cong P/M$.\findemo
\section{Simple biset functors and bilinear forms}
Let $\F$ be any field. Recall (see~\cite{boustathe}) that, given a finite group $H$, we defined, for any finite group $G$
$$\F\sur{B}(G,H)=\F B(G,H)/\sum_{|K|<|H|}\F B(G,K)\circ \F B(K,H)\mvirg$$
and that the correspondence $G\mapsto \F \sur{B}(G,H)$ is a quotient biset functor of the Yoneda functor $G\mapsto \F B(G,H)$ at the group $H$.\par
When $V$ is a $\F\Out(H)$-module, we defined an $\F$-valued bilinear form $\scalv{\;}{\,}_{V,G}$ on $\F \sur{B}(G,H)$ by
$$\forall \alpha,\beta\in\F \sur{B}(G,H),\;\scalv{\alpha}{\beta}_{V,G}=\chi_V\big(\pi_H(\hat{\alpha}\op\circ\hat{\beta})\big)\mvirg$$
where $\hat{\alpha}, \hat{\beta}$ are elements of $\F B(G,H)$ lifting $\alpha, \beta\in\F \sur{B}(G,H)$, respectively, where $\pi_H:\F B(H,H)\to \F \sur{B}(H,H)\cong \F \Out(H)$ is the projection map, and $\chi_V$ is the character of $V$, i.e. the trace function $\End_\F(V)\to \F$. The main property of these constructions is that
$$\F \sur{B}(G,H)/{\rm Rad}\scalv{\;}{\,}_{V,G}\cong S_{H,V}(G)^{\dim_\F V}\mpoint$$
Moreover, if $L$ is a finite group, then for any $\gamma\in B(L,G)$, any $\alpha\in\sur{B}(G,H),$ and any $\beta\in\sur{B}(L,H)$,
\begin{equation}\label{adjoint}
\scalv{\gamma(\alpha)}{\beta}_{V,L}=\scalv{\alpha}{\gamma\op(\beta)}_{V,G}\mpoint
\end{equation}
\vspace{2ex}
\pagebreak[3]
\npar Suppose from now on that $\F$ is a field of characteristic 0. Observe that $\tilde{e}_K^G=(\tilde{e}_K^G)\op$ for any subgroup $K$ of a finite group $G$. By~(\ref{adjoint}), this implies that the decomposition 
$$\F\sur{B}(G,H)=\bigoplus_{K\in[s_G]}\tilde{e}_K^G\F\sur{B}(G,H)\mvirg$$
where $[s_G]$ is a set of representatives of conjugacy classes of subgroups of $G$, is an orthogonal decomposition with respect to the form $\scalv{\;}{\,}_{V,G}$. Moreover
$$\tilde{e}_G^GS_{H,V}(G)^{\dim_\F V}\cong \tilde{e}_G^G\F\sur{B}(G,H)/{\rm Rad}\scalv{\;}{\,}_{V,G}\mvirg$$
and the isomorphism
$$
\F\sur{B}(G,H)\cong\bigoplus_{K\in[s_G]}\big(\tilde{e}_K^K\F\sur{B}(K,H)\big)^{N_G(K)}
$$
given by Proposition~6.5.5 of~\cite{bisetfunctorsMSC} induces an isomorphism
\begin{equation}
\label{direct sum decomposition}S_{H,V}(G)^{\dim_\F V}\cong \bigoplus_{K\in[s_G]}\big(\tilde{e}_K^K\F\sur{B}(K,H)/{\rm Rad}\scalv{\;}{\,}_{V,K}\big)^{N_G(K)}\mpoint
\end{equation}
Now $\sur{B}(K,H)$ is generated by the images of the elements $(K\times H)/L$, where $L$ is a subgroup of $K\times H$. If this image is non-zero, then $L$ is of the form $L=\{\big(x,s(x)\big)\mid x\in X\}$, where $X$ is a subgroup of $K$ and $s:X\twoheadrightarrow H$ is a surjective group homomorphism. \par
The $(K,G)$-biset $U=(K\times H)/L$ factors as $U=\Ind_X^K\circ V$, for a suitable $(X,H)$-biset $V$ (\cite{bisetfunctorsMSC}, Lemma~2.3.26), and by \cite{bisetfunctorsMSC}, Corollary 2.5.12 
$$\tilde{e}_K^K\circ \Ind_X^K=\Ind_X^K\circ \widetilde{\Res_X^Ke_K^K}\mpoint$$
Now $\Res_X^Ke_K^K=0$ if $X$ is a proper subgroup of $K$. It follows that $\tilde{e}_K^K\F\sur{B}(K,H)$ is generated by the images $\sur{u}_s$ of the elements
$$u_s=\tilde{e}_K^K\times_K(K\times H)/\Delta_s^\diamond (K)\mvirg$$
where $\Delta_s^\diamond (K)=\{\big(x,s(x)\big)\mid x\in K\}$, for a surjective group homomorphism $s:K\twoheadrightarrow H$.\par
Let $\varpi_H:\Aut(H)\to\Out(H)$ denote the projection map. Then:
\pagebreak[3]
\begin{mth}{Proposition} \label{the formula}Let $s,t:K\twoheadrightarrow H$ be two surjective group homomorphisms. Let $M=\Ker\,s$, and $N=\Ker\,t$. Then
\begin{eqnarray*}
\scalv{\sur{u}_s}{\sur{u}_t}_{V,K}\!&\!=\!&\!m_{K,M\cap N}\,\fracb{\mu_{\normal K}(M\cap N,M)}{|M:M\cap N|}\,\chi_V\big(\sum_{Y\in\sur{\mathcal{K}}(K,M,N)}\varpi_H([s,Y,t])\big)\\
&\!=\!&\!m_{K,M\cap N}\,\fracb{\mu_{\normal K}(M\cap N,M)}{|M:M\cap N|}\,\chi_V\big(\sumb{\theta\in \Aut(H)}{\Delta_\theta(H)\leq (s\times t)(K)}\varpi_H(\theta)\big)\mvirg
\end{eqnarray*}
where $\mu_{\normal K}$ is the Möbius function of the poset of normal subgroups of $K$, and
$$\sur{\mathcal{K}}(K,M,N)=\{Y\leq K\mid YN=YM=K,\;Y\cap N=Y\cap M=M\cap N\}$$ 
is the set of subgroups $Y$ of $K$, containing $M\cap N$, such that $Y/(M\cap N)$ is a common complement of $M/(M\cap N)$ and $N/(M\cap N)$ in $K/(M\cap N)$. Moreover for $Y\in\sur{\mathcal{K}}(K,M,N)$, the symbol $[s,Y,t]$ denotes the automorphism of $H$ defined by $[s,Y,t]\big(t(y)\big)=s(y)$, $\forall y\in Y$. 
\end{mth}
\pf By definition, and since $\tilde{e}_K^K$ is an idempotent
\begin{eqnarray*}
\scalv{\sur{u}_s}{\sur{u}_t}_{V,K}&=&\chi_V\big(\pi_H(u_s\op\circ u_t)\big)\\
&=&\chi_V\Big(\pi_H\big((H\times K)/\Delta_s(K)\times_K\tilde{e}_K^K\times_K(K\times H)/\Delta_t^\diamond(K)\big)\Big)\mvirg
\end{eqnarray*}
where $\Delta_s(K)=\{\big(s(x),x\big)\mid x\in K\}$. Set
$$a_{s,t}=(H\times K)/\Delta_s(K)\times_K\tilde{e}_K^K\times_K(K\times H)/\Delta_t^\diamond(K)\mpoint$$
Then
$$a_{s,t}=\frac{1}{|K|}\sum_{L\leq K}|L|\mu(L,K)(H\times H)/\Delta_{s,t}(L)\mvirg$$
where $\Delta_{s,t}(L)=\{\big(s(l),t(l)\big)\mid l\in L\}$.\par
This subgroup of $H\times H$ is equal to $\Delta_\theta(H)$, for some automorphism $\theta$ of~$H$, if and only if 
\begin{equation}\label{conditions}L\cap M=L\cap N\;\;\hbox{and}\;\;LM=K=LN\mvirg\end{equation}
where $M=\Ker\,s$ and $N=\Ker\,t$. In this case the automorphism $\theta$ is defined by $\theta\big(t(l)\big)=s(l)$, for any $l\in L$. The two conditions~\ref{conditions} and the automorphism $\theta$ remain unchanged when $L$ is replaced by $L(M\cap N)$. Moreover the conditions~\ref{conditions} are equivalent to saying that the group $Y=L(M\cap N)$ is in $\sur{\mathcal{K}}(K,M,N)$, and in this case $\theta=[s,Y,t]$. Conversely, fix some $Y\in\sur{\mathcal{K}}(K,M,N)$, and consider all the subgroups $L$ of $K$ such that $L(M\cap N)=Y$. Recall that
$$\sumb{L\leq Y}{L(M\cap N)=Y}|L|\mu(L,K)=m_{K,M\cap N}|Y|\mu(Y,K)\mpoint$$
This gives:
$$\pi_H(a_{s,t})=m_{K,M\cap N}\sum_{Y\in \sur{\mathcal{K}}(K,M,N)}\frac{|Y|}{|K|}\mu(Y,K)\varpi([s,Y,t])\mpoint$$
Now if $Y/(M\cap N)$ is a complement of $M/(M\cap N)$ in $K/(M\cap N)$, and if $K/M\cong H$, it follows that $|Y|=|M\cap N||H|$. Moreover the poset $]Y,K[$ is isomorphic to the poset $]M\cap N,M|^Y$. But since $M$ and $N$ are normal subgroups of $K$, the commutator group $[M,N]$ is contained in $M\cap N$. It follows that $]M\cap N,M|^Y=]M\cap N,M|^{YN}=]M\cap N,M|^K$, and that $\mu(Y,K)=\mu_{\normal K}(M\cap N,M)$. This completes the proof of the first equality of the proposition. The second one follows from the observation that the correspondences
$$Y\mapsto [s,Y,t]\;\;\hbox{and}\;\;\theta\mapsto \{k\in K\mid \theta\big(t(k)\big)=s(k)\}$$
are mutual inverse bijections between $\sur{\mathcal{K}}(K,M,N)$ and the set of automorphisms $\theta$ of $H$ such that $\Delta_\theta(H)\leq (s\times t)(K)$ (see Section 8.3 of~\cite{doublact}).\findemo
\begin{mth}{Corollary} \label{reduction}\begin{enumerate}
\item If $\tilde{e}_K^KS_{H,V}(K)\neq\zero$, the group $\beta(K)$ is isomorphic to $\beta(L)$, where $L$ is a subgroup of $H\times H$ with the following properties:
\begin{enumerate}
\item $p_1(L)=p_2(L)=H$. 
\item $k_1(L)$ and $k_2(L)$ are direct products of minimal normal subgroups of $H$.
\item There exist an automorphism $\theta$ of $H$ such that $\theta\big(k_2(L)\big)=k_1(L)$.
\end{enumerate}
\item In particular, if $H$ is a $p$-group for some prime $p$, then $K$ is $p$-elementary. If $K=P\times C$, where $P$ is a $p$-group and $C$ is a cyclic $p'$-group, then
$$\tilde{e}_K^KS_{H,V}(K)\cong \tilde{e}_P^PS_{H,V}(P)\mvirg$$
and this isomorphism is compatible with the action of $\Aut(K)$.
\end{enumerate}
\end{mth}
\pf Indeed, if $\tilde{e}_K^KS_{H,V}(K)\neq\zero$, then the bilinear form $\scalv{\;}{\,}_{V,K}$ is not identically zero on $\tilde{e}_K^K\F\sur{B}(K,H)$. It follows that there exist surjective group homomorphisms $s,t:K\twoheadrightarrow H$ such that $\scalv{\sur{u}_s}{\sur{u}_t}_{V,K}\neq 0$. \par
Then $m_{K,M\cap N}\neq 0$, $\mu_{\normal K}(M\cap N,M)\neq 0$, and $\sur{\mathcal{K}}(K,M,N)\neq \emptyset$, where $M=\Ker\,s$ and $N=\Ker\,t$. Hence $\beta(K)\cong\beta\big(K/(M\cap N)\big)$. Now $K/(M\cap N)$ is isomorphic to $L=(s\times t)(K)$, which is a subgroup of $H\times H$ such that $p_1(L)=p_2(L)=H$. Moreover $k_1(L)=s(\Ker\,t)\cong N/(M\cap N)$ and $k_2(L)=t(\Ker\,s)\cong M/(M\cap N)$. Then $\mu_{\normal K}(M\cap N,M)=\mu_{\normal H}\big(\un,t(\Ker\,s)\big)$, and this is non zero if and only if the lattice $[\un,t(\Ker\,s)]^H$ of normal subgroups of $H$ contained in $t(\Ker\,s)$ is complemented, i.e. if $t(\Ker\,s)$ is a direct product of minimal normal subgroups of $H$. Finally, let $Y\in\sur{\mathcal{K}}(K,M,N)$ and $\theta=[s,Y,t]$. If $u\in k_2(L)=t(\Ker\,s)$, then there exist $v\in \Ker\,s)$ and $y\in Y$ such that $u=t(v)=t(y)$. Then $v^{-1}y\in \Ker\,t$, and $\theta(u)=s(y)=s(v^{-1}y)\in s(\Ker\,t)=k_1(L)$. In other words $\theta\big(k_2(L))=k_1(L)$, which completes the proof of Assertion~1.\par 
The first part of Assertion~2 follows from Assertion~2 of Lemma~\ref{p-elementary}. Now if $H$ is a $p$-group, if $K=P\times C$, where $P$ is a $p$-group and $C$ is a $p'$-group, and if $s:K\twoheadrightarrow H$ is a surjective group homomorphism, then $C\leq \Ker\,s$. In other words, there is a surjective homomorphism $\sur{s}:P\twoheadrightarrow H$ such that $s=\sur{s}\circ \pi$, where $\pi :K\to P$ is the projection map. Moreover $\Ker\,s=\Ker\,\sur{s}\times C$.\par
So with the notation of Proposition~\ref{the formula}, $M=\sur{M}\times C$, where $\sur{M}=\Ker\,\sur{s}$. Similarly $N=\sur{N}\times C$, where $N=\Ker\,\sur{t}$, and $\sur{t}:P\twoheadrightarrow H$ is such that $t=\sur{t}\circ \pi$. Clearly $|M:M\cap N|=|\sur{M}:\sur{M}\cap\sur{N}|$. Moreover, one checks easily that
$$m_{K,M\cap N}=m_{P\times C,(\sur{M}\cap\sur{N})\times C}=m_{P,\sur{M}\cap\sur{N}}m_{C,C}=m_{P,\sur{M}\cap\sur{N}}\fracb{\phi(|C|)}{|C|}\mvirg$$
since $P$ and $C$ have coprime orders, and $C$ is cyclic. \par
Also
$$\mu_{\normal K}(M\cap N,M)=\mu_{\normal P}(\sur{M}\cap\sur{N},\sur{M})\mpoint$$
Finally, the maps $Q\mapsto Q\times\un$ and $Y\mapsto (Y\cap P)$ induce inverse bijections from $\sur{\mathcal{K}}(P,\sur{M},\sur{N})$ to $\sur{\mathcal{K}}(K,M,N)$, and for any $Y\in \sur{\mathcal{K}}(K,M,N)$,
$$[s,Y,t]=[\,\sur{s},Y\cap P,\sur{t}\,]\mpoint$$
It follows that the matrix of the form $\scalv{\;}{\,}_{V,K}$ on $\tilde{e}_K^K\F\sur{B}(K,H)$ is equal to the matrix of the form $\scalv{\;}{\,}_{V,P}$ on $\tilde{e}_P^P\F\sur{B}(P,H)$, multiplied by the non-zero scalar $\frac{\phi(|C|)}{|C|}$. Hence the two forms define isomorphic quadratic spaces. As all the above bijections are obviously compatible with the action of $\Aut(K)$ and the canonical group homomorphism $\Aut(K)\to\Aut(P)$, the induced isomorphism
$$\tilde{e}_K^KS_{H,V}(K)\cong \tilde{e}_P^PS_{H,V}(P)$$
is compatible with the action of $\Aut(K)$.\findemo
\begin{mth}{Notation} \label{divers} Let $H$ and $P$ be finite $p$-groups. 
\begin{enumerate}
\item Let $\mathcal{Q}_{H}(P)$ denote the $\F$-vector space with basis the set
$$\Sigma_{H}(P)=\{s\mid s:P\twoheadrightarrow H\}$$
of surjective group homomorphisms from $P$ to $H$, endowed with the $\F$-valued bilinear form $\scalv{\;}{\,}_{V,P}$ defined as follows: For $s,t\in\Sigma_{H}(P)$, set $M=\Ker\,s$ and $N=\Ker\,t$. If $M\cap\Phi(P)\neq N\cap\Phi(P)$, set $\scalv{s}{t}_{V,P}=0$. And if $M\cap\Phi(P)=N\cap\Phi(P)$, then the groups $M/(M\cap N)$ and $N/(M\cap N)$ are central elementary abelian subgroups of the same rank of $P/M\cap N$. In this case, set
$$\scalv{s}{t}_{V,P}=m_{P,M\cap N}\fracb{\mu(M\cap N,M)}{|M:M\cap N|}\chi_V\big(\sum_{Y\in\sur{\mathcal{K}}(P,M,N)}\varpi_H([s,Y,t])\big)\mpoint$$
\item Let $\mathcal{Q}_{H}^\sharp(P)$ be the subspace of $\mathcal{Q}_{H}(P)$ with basis the subset
$$\Sigma_{H}^\sharp(P)=\{s\mid s:P\twoheadrightarrow H,\;\;\Ker\,s\cap\Phi(P)=\un\}$$
of $\Sigma_H(P)$.
\item Set 
$$\mathcal{N}_H(P)=\{N\mid N\normal P,\;\;N\cap\Phi(P)=\un\}\mpoint$$
\item Denote by $\mathcal{E}_H(P)$ the set of normal subgroups $R$ of $P$, contained in $\Phi(P)$, and such that $P/R\cong E\times H$, for some elementary abelian $p$-group $E$.
\end{enumerate}
\end{mth}
\begin{mth}{Proposition} \label{blocs}Let $H$ be a $p$-group, and $K$ be a $p$-elementary group. Set $P=O_p(K)$. Then:
\begin{enumerate}
\item There is an isomorphism of $\F\Aut(K)$-modules
$$\tilde{e}_K^KS_{H,V}(K)\cong\mathcal{Q}_{H}(P)/{\rm Rad}\scalv{\;}{\,}_{V,P}\mpoint$$
\item Let $\Gamma$ be a finite group acting on the group $K$. Then $\Gamma$ acts on the set $\mathcal{E}_H(P)$, and there is an isomorphism of $\F\Gamma$-modules
$$\tilde{e}_K^KS_{H,V}(K)\cong\dirsum{R\in[\Gamma\dom\mathcal{E}_H(P)]}\Ind_{\Gamma_R}^{\Gamma}\big(\mathcal{Q}_{H}^\sharp(P/R)/{\rm Rad}\scalv{\;}{\,}_{V,P/R}\big)\mvirg$$
where  $[\Gamma\dom\mathcal{E}_H(P)]$ is a set of representatives of $\Gamma$-orbits on $\mathcal{E}_H(P)$, and $\Gamma_R$ denotes the stabilizer of~$R$ in $\Gamma$.
\end{enumerate}
\end{mth}
\pf The map $s\in\Sigma_H(P)\mapsto \sur{u}_s\in \tilde{e}_P^P\F\sur{B}(P,H)$ induces a surjective linear map $\mathcal{Q}_H(P)\to \tilde{e}_P^P\F\sur{B}(P,H)$. 
Let $s,t\in\Sigma_H(P)$, and set $M=\Ker\,s$ and $N=\Ker\,t$. Then $|M|=|N|$. It follows from Lemma~\ref{complements} that $\sur{\mathcal{K}}(P,M,N)\neq\emptyset$ if and only if $M/(M\cap N)$ and $N/(M\cap N)$ are central elementary abelian subgroups of $P/(M\cap N)$, which intersect trivially the Frattini subgroup of $P/(M\cap N)$. But
$$\Phi\big(P/(M\cap N)\big)=\Phi(P)(M\cap N)/(M\cap N)\mpoint$$
Hence $M/(M\cap N)\cap \Phi\big(P/(M\cap N)\big)=\big(M\cap \Phi(P)\big)(M\cap N)/(M\cap N)$. This group is trivial if and only if $M\cap\Phi(P)\leq M\cap N$, i.e. if $M\cap\Phi(P)\leq N\cap \Phi(P)$. Hence $\sur{\mathcal{K}}(P,M,N)\neq\emptyset$ if and only if $M\cap\Phi(P)=N\cap\Phi(P)$.\par
 This holds in particular if $\scalv{\sur{u}_s}{\sur{u}_t}\neq 0$. In this case, by Proposition~\ref{the formula} 
$$ \scalv{\sur{u}_s}{\sur{u}_t}_{V,P}=m_{P,M\cap N}\,\fracb{\mu_{\normal P}(M\cap N,M)}{|M:M\cap N|}\,\chi_V\big(\sum_{Y\in\sur{\mathcal{K}}(P,M,N)}\varpi_H([s,Y,t])\big)\mpoint$$
But $\mu_{\normal P}(M\cap N,M)=\mu(M\cap N,M)$, since $P/(M\cap N)$ centralizes both $M/(M\cap N)$ and $N/(M\cap N)$. Hence
$$\scalv{\sur{u}_s}{\sur{u}_t}_{V,P}=\scalv{s}{t}_{V,P}$$
in this case, and Assertion~1 follows.\par
Since $\scalv{s}{t}_{V,P}=0$ if $M\cap\Phi(P)\neq N\cap\Phi(P)$, the quadratic space $\mathcal{Q}=\big(\mathcal{Q}_H(P),\scal{\;}{\;}_{V,P}\big)$ splits as the orthogonal sum of the subspaces $\mathcal{Q}_R$ generated by the elements $s\in\Sigma_H(P)$ such that $\Ker\,s\cap\Phi(P)=R$. These subspaces are permuted by the action of $\Aut(K)$, and the space $\mathcal{Q}_R$ is invariant by $\Aut(K)_R$. \par
Let $\pi_R:P\to P/R$ be the canonical projection. The map
$$\theta_R:\sur{s}\in\Sigma_H^\sharp(P/R)\mapsto \sur{s}\circ\pi_R$$
is a bijection from $\Sigma_H^\sharp(P/R)$ to the set $\{s\in\Sigma_H(P)\mid \Ker\,s\cap\Phi(P)=R\}$, and the map $Y\mapsto Y/R$ is a bijection from $\sur{\mathcal{K}}(P,M,N)$ to $\sur{\mathcal{K}}(P/R,M/R,N/R)$, such that
$$[\sur{s},Y/R,\sur{\,t}\,]=[\theta_R(\sur{s}),Y,\theta_R(\sur{\,t}\,)]\mvirg$$
for any $\sur{s},\sur{\,t}\,\in\Sigma_H^\sharp(P/R)$.\par
Moreover, if $M\cap\Phi(P)=N\cap \Phi(P)=R$, then
$$m_{P,M\cap N}=m_{P,R}m_{P/R,(M\cap N)/R}=m_{P/R,(M\cap N)/R}\mvirg$$
as $R\leq \Phi(P)$. Also
$$M/(M\cap N)\cong (M/R)/\big((M/R)\cap(N/R)\big)\mpoint$$
It follows that
$$\forall \sur{s},\sur{\,t}\,\in \Sigma_H^\sharp(P/R),\;\;\scalv{\theta_R(\sur{s})}{\theta_R(\sur{\,t}\,)}_{V,P}=\scalv{\sur{s}}{\sur{\,t}\,}_{V,P/R}\mpoint$$
Hence there is an isomorphism
 $$\mathcal{Q}_R/{\rm Rad}\scalv{\;}{\,}_{V,P}\cong \mathcal{Q}_H^\sharp(P/R)/{\rm Rad}\scalv{\;}{\,}_{V,P/R}\mvirg$$
of $\F\Gamma_R$-modules. \par
To complete the proof of Assertion~2, it remains to observe that if $Q$ is a $p$-group, the set $\Sigma_H^\sharp(Q)$ is non-empty if and only if the group $Q$ is isomorphic to $E\times H$, for some elementary abelian $p$-group $E$: Indeed, if $Q=E\times H$, where $E$ is elementary abelian, then $\Phi(Q)=\un\times \Phi(H)$, and the projection map $s:Q\to H=Q/E$ is an element of $\Sigma_H^\sharp(Q)$. Conversely, if $s\in\Sigma_H^\sharp(Q)$, then $E=\Ker\,s$ is an elementary abelian central subgroup of $Q$, which admits a complement $L$ in $Q$, by Lemma~\ref{complements}. Thus $Q=E\times L$, and $L\cong Q/E\cong H$. Hence $Q\cong E\times H$.
 \findemo
 \begin{mth}{Theorem} Let $G$ be a finite group, let $H$ be a finite $p$-group, and let $V$ be a simple $\F\Out(H)$-module. Then
 $$S_{H,V}(G)^{\dim_\F V}\cong\dirsum{(K,R)}\Big(\mathcal{Q}_H^\sharp(K_p/R)/{\rm Rad}\scalv{\;}{\,}_{V,K_p/R}\Big)^{N_G(K,R)}\mpoint$$
 where $(K,R)$ runs through a set of $G$-conjugacy classes of pairs consisting of a $p$-elementary subgroup $K$ of $G$, and a $p$-subgroup $R$ in $\mathcal{E}_H(K_p)$, where $K_p=O_p(K)$, and $N_G(K,R)=N_G(K)\cap N_G(R)$.
 \end{mth}
\pf This follows from Equation~\ref{direct sum decomposition}, Corollary~\ref{reduction} and Proposition~\ref{blocs}.\findemo
\section{Proof of the theorem}
This section is devoted to the proof of the following theorem, announced in the introduction:
\begin{mth}{Theorem} \label{the theorem}Let $\F$ be a field of characteristic 0, let $p$ be a prime number, and $H$ be a finite $p$-group. Let moreover $G$ be a finite group.
\begin{enumerate}
\item If $H=\un$, the dimension of $S_{H,\F}(G)$ is equal to the number of conjugacy classes of cyclic subgroups of $G$.
\item If $H\cong C_p\times C_p$, the dimension of $S_{H,\F}(G)$ is equal to the number of conjugacy classes of non-cyclic $p$-elementary subgroups of $G$.
\item If $H$ is any other finite $p$-group, the dimension of $S_{H,\F}(G)$ is equal to the number of conjugacy classes of sections $(T,S)$ of $G$ such that $T/S\cong H$ and $T$ is $p$-elementary.
\end{enumerate}
\end{mth}
\pf {\bf Step 1:} Let $K$ be a subgroup of $G$. Then $\tilde{e}_K^KS_{H,\F}(K)=\zero$, by Corollary~\ref{reduction}, unless $K\cong P\times C$, where $P$ is a $p$-group and $C$ is a cyclic $p'$-group, and in this case $\tilde{e}_K^KS_{H,\F}(K)\cong \tilde{e}_P^PS_{H,\F}(P)$ as $\F\Aut(K)$-modules.\par
By Proposition~\ref{blocs}, there is an isomorphism of $\F\Aut(P)$-modules
$$\tilde{e}_P^PS_{H,\F}(P)\cong\dirsum{R}\Ind_{\Aut(P)_R}^{\Aut(P)}\big(\mathcal{Q}_{H}^\sharp(P/R)/{\rm Rad}\scalv{\;}{\,}_{\F,P/R}\big)\mvirg$$
where $R$ runs through a set of representatives of $\Aut(P)$-orbits of normal subgroups of $P$ contained in $\Phi(P)$, such that $P/R\cong E\times H$, for some elementary abelian $p$-group $E$. So the computation of $\tilde{e}_P^PS_{H,\F}(P)$ comes down to the computation of the $\F\Aut(Q)$-module
$$\mathcal{V}_H(Q)=\mathcal{Q}_{H}^\sharp(Q)/{\rm Rad}\scalv{\;}{\,}_{\F,Q}\mvirg$$
for a $p$-group $Q=P/R$ of the form $E\times H$, where $R$ is some normal subgroup of $P$ contained in $\Phi(P)$.
Recall that $\mathcal{Q}_{H}^\sharp(Q)$ is the $\F$-vector space with basis 
$$\Sigma_H^\sharp(Q)=\{s\mid s:Q\twoheadrightarrow H,\;\;\Ker\,s\cap\Phi(Q)=\un\}\mvirg$$
and that the bilinear form $\scalv{\;}{\,}_{\F,Q}$ is defined for $s,t\in \Sigma_H^\sharp(Q)$ by
$$\scalv{s}{t}_{\F,Q}=m_{Q,M\cap N}\fracb{\mu(M\cap N,M)}{|M:M\cap N|}|\sur{\mathcal{K}}(Q,M,N)|\mvirg$$
where $M=\Ker\,s$ and $N=\Ker\,t$.\par
This shows that $\scalv{s}{t}_{\F,Q}$ depends only on $M$ and $N$. It follows that $\mathcal{Q}_H^\sharp(Q)/{\rm Rad}\scalv{\;}{\,}_{\F,Q}$ is also isomorphic to the quotient of the $\F$-vector space with basis the set $\mathcal{N}_H(Q)=\{N\normal Q\mid P/N\cong H,\;\;N\cap\Phi(Q)=\un\}$ introduced in Notation~\ref{divers}, by the radical of the bilinear form $\ncalv{\;}{\,}_{\F,Q}$ defined by
 $$\ncalv{M}{N}_{\F,Q}=m_{Q,M\cap N}\,\fracb{\mu(M\cap N,M)}{|M:M\cap N|}\,|\sur{\mathcal{K}}(Q,M,N)|\mvirg$$
for $M, N\in\mathcal{N}_H(Q)$ .\spn
{\bf Step 2:} Now Assertion~1 is well known (see e.g. Proposition 4.4.8~\cite{bisetfunctorsMSC}), but it can also be recovered from the argument of Step~1: Indeed, if $H=\un$, there is a unique normal subgroup $N$ of $Q$ such that $Q/N\cong H$, namely $Q$ itself. If moreover $N\cap \Phi(Q)=1$, then $\Phi(Q)=\un$, and $Q$ is elementary abelian. But as $Q=P/R$, for some $R\leq\Phi(P)$, it follows that $R=\Phi(P)$, and $Q=P/\Phi(P)$. Moreover $\mathcal{N}_H(Q)=\{Q\}$, and
$$\ncalv{Q}{Q}_{\F,Q}=m_{Q,Q}\mvirg$$
which is equal to 0 if $Q$ is non-cyclic, and to $1-1/p$ otherwise. Hence $\mathcal{V}_H(Q)=\zero$ if $Q$ is non-cyclic, and $\mathcal{V}_H(Q)$ is one dimensional if $Q=P/\Phi(P)$ is cyclic, i.e. if $P$ is cyclic. But $P=O_p(K)$ for some $p$-elementary subgroup $K$ of $G$. Hence $P$ is cyclic if and only if $K$ itself is cyclic, and this leads to Assertion~1.\sp
We can now assume that $H$ is a non-trivial $p$-group, of order $p^h$, and make a series of observations:
\begin{itemize}
\item Let $P$ be a $p$-group, and $Q$ be a normal subgroup of $P$. Example~5.2.3 of~\cite{bisetfunctorsMSC} shows that $m_{P,Q}=m_{P,\Phi(P)Q}$. By Proposition~5.3.1 of~\cite{bisetfunctorsMSC}, it follows that
\begin{equation}\label{mPQ}
m_{P,Q}=m_{P,\Phi(P)}m_{P/\Phi(P),Q\Phi(P)\Phi(P)}=m_{E,F}\mvirg
\end{equation}
where $E$ is the elementary abelian $p$-group $P/\Phi(P)$, and $F$ its subgroup $Q\Phi(P)/\Phi(P)$. If $E$ has rank $n\geq 2$ and $F$ has rank $k$, then
\begin{equation}\label{mEF}m_{E,F}=(1-p^{n-2})(1-p^{n-3})\cdots(1-p^{n-k-1})\mvirg\end{equation}
(which is equal to 0 if $k\geq n-1$, and non-zero otherwise). This follows from an easy induction argument on $k$, using Proposition~5.3.1 of~\cite{bisetfunctorsMSC}, and starting with the case $k=1$, which is a special case of Equation~\ref{minimal normal}.\par
If $E$ has rank 1 and $F=\un$, then $m_{E,F}=1$. In this case $m_{E,E}=1-\frac{1}{p}$. This is the only case where $m_{E,F}$ is not an integer.
\item  Let $M,N\in\mathcal{N}_H(Q)$. Then in particular $M$ and $N$ have the same order. Recall that $m_{Q,M\cap N}$ is non-zero if and only if $\beta(Q)\cong\beta\big(Q/(M\cap N))$. So either $Q$ and $Q/(M\cap N)$ are both cyclic, or they are both non-cyclic. Equivalently, either $Q$ is cyclic, or $Q/(M\cap N)$ is non-cyclic. If $H$ is non-cyclic, then $Q/(M\cap N)$ is non-cyclic, as it maps surjectively on $Q/M\cong H$. \par
So if $m_{Q,M\cap N}=0$, then $H$ is cyclic, $Q$ is non-cyclic, and $Q/(M\cap N)$ is cyclic. But then $M/(M\cap N)=N/(M\cap N)$, since the cyclic group $Q/(M\cap N)$ admits a unique subgroup of a given order. Thus $M=N$. Conversely, if $H$ is cyclic, if $Q$ is non-cyclic, and if $M=N$, then $m_{Q,M\cap N}=0$ since $Q/(M\cap N)\cong H$ is cyclic.
\item If $M,N\in\mathcal{N}_H(Q)$, then the subgroups $\sur{M}=M/(M\cap N)$ and $\sur{N}=N/(M\cap N)$  are central elementary abelian subgroups of the same order of $\sur{Q}=Q/(M\cap N)$. 
If $\sur{M}$ has rank $m$, then
$$\mu(M\cap N,M)=(-1)^mp^{\binom{m}{2}}\mpoint$$
Now  by Lemma~\ref{complements}, the product
\begin{equation}\label{alphaMN}
\alpha_{M,N}=\mu(M\cap N,M)|\sur{K}(Q,M,N)|
\end{equation}
is equal to
\begin{eqnarray*}
\alpha_{M,N}&=&(-1)^m(p^s-1)(p^{s-1}-1)\cdots(p-1)p^{\binom{m}{2}+\binom{s}{2}+s(m-s)+m(\gamma-s)}\\
&=&(-1)^{m+s}(1-p^s)(1-p^{s-1})\cdots(1-p)p^{\binom{m}{2}+\binom{s}{2}+s(m-s)+m(\gamma-s)}
\mvirg
\end{eqnarray*}
where $\gamma$ is the rank of $H/\Phi(H)$, and $s$ is the rank of $\sur{M}/\big(\sur{M}\cap \sur{N}\Phi(\sur{Q})\big)\cong M/\big(M\cap N\Phi(Q)\big)$.
\end{itemize}
{\bf Step 3:} Finally $\ncalv{M}{N}_{\F,Q}=0$ if and only if $m_{Q,M\cap N}=0$, i.e. $H$ is cyclic, $M=N$, and $Q$ is not cyclic. In all other cases, the groups $\sur{M}=M/(M\cap N)$ and $\sur{N}=N/(M\cap N)$ are elementary abelian, and central in $\sur{Q}=Q/(M\cap N)$. Moreover $\sur{M}\cap\sur{N}=\un$. Let $m$ be the rank of~$\sur{M}$, let $s$ denote the rank of~$\sur{M}/\big(\sur{M}\cap \sur{N}\Phi(\sur{Q})\big)$, and let $\gamma$ denote the rank of $H/\Phi(H)$. Then
\begin{eqnarray*}\ncalv{M}{N}_{\F,Q}\!\!&\!\!=\!\!&\!\!m_{Q,M\cap N}\fracb{\alpha_{M,N}}{|M:M\cap N|}\\
\!\!&\!\!=\!\!&\!\!m_{Q,M\cap N}\,\fracb{\alpha_{M,N}}{|\sur{M}|}\\
\!\!&\!\!=\!\!&\!\!(-1)^m m_{Q,M\cap N}(p^s\!-\!1)(p^{s-1}\!-\!1)\cdots(p\!-\!1)p^{\binom{m}{2}-m+\binom{s}{2}+s(m-s)+m(\gamma-s)}\mvirg\\
\end{eqnarray*}
i.e. finally
\begin{equation}\label{scalar}
\ncalv{M}{N}_{\F,Q}=(-1)^m m_{Q,M\cap N}(p^s\!-\!1)(p^{s-1}\!-\!1)\cdots(p\!-\!1)p^{\frac{1}{2}{(m-s)(m+s+1)}+m(\gamma-2)}\mpoint
\end{equation}
Let $n$ denote the rank of $Q/\Phi(Q)$. By Equation~\ref{mEF}
$$m_{Q,M\cap N}=(1-p^{n-2})(1-p^{n-3})\cdots(1-p^{n-k-1})\mvirg$$
where $k$ is the rank of $(M\cap N)\Phi(Q)/\Phi(Q)\cong M\cap N$. Since $Q/M\Phi(Q)\cong H/\Phi(H)$ has rank $\gamma$, it follows that $M\Phi(Q)/\Phi(Q)\cong M$ has rank $n-\gamma$. Since $M/(M\cap N)$ has rank $m$, it follows that $k=n-m-\gamma$. Thus
$$m_{Q,M\cap N}=(1-p^{n-2})(1-p^{n-3})\cdots(1-p^{m+\gamma-1})\mpoint$$
It follows that
\begin{equation}\label{scalar2}\ncalv{M}{N}_{\F,Q}=A_{M,N}(-1)^{m+s}p^{\frac{1}{2}{(m-s)(m+s+1)}+m(\gamma-2)}\mvirg
\end{equation}
where 
$$A_{M,N}=(1-p^{n-2})(1-p^{n-3})\cdots(1-p^{m+\gamma-1})(1-p^s)(1-p^{s-1})\cdots(1-p)$$ 
is an integer congruent to 1 modulo $p$.\spn
\noindent{\bf Step 4:} Assume first that $H$ is non-cyclic, i.e. that $\gamma\geq 2$. In this case $\ncalv{M}{N}_{\F,Q}$ is non-zero.
If $M=N$, then $m=s=0$, and $\ncalv{M}{N}_{\F,P}=A_{M,M}$ is congruent to 1 modulo $p$. And if $M\neq N$, then $m\geq 1$. As $\gamma\geq 2$ and $m\geq s$, the exponent
$$\frac{1}{2}{(m-s)(m+s+1)}+m(\gamma-2)$$
of $p$ in the right hand side of~\ref{scalar2} is non-negative. It is equal to 0 if and only if $m=s$ and $\gamma=2$. In this case $\sur{M}\cap \sur{N}\Phi(\sur{Q})=\un$, so $\sur{M}$ maps into $\sur{Q}/\sur{N}\Phi(\sur{Q})\cong H/\Phi(H)$, which has rank $\gamma=2$. It follows that $m\leq 2$. \par
If $m=2$, then $\sur{M}\sur{N}\Phi(\sur{Q})=\sur{Q}$, thus $\sur{M}\sur{N}=\sur{Q}$, and $H\cong \sur{Q}/\sur{N}\cong\sur{M}$ (since $\sur{M}\cap\sur{N}=\un$), so $H$ is elementary abelian of rank 2.\par
If $m=1$, then as $\sur{M}\cong C_p$ maps into $\sur{Q}/\sur{N}\Phi(\sur{Q})\cong C_p\times C_p$, the group $\sur{Q}/\big(\sur{M}\sur{N}\Phi(\sur{Q})\big)$ is cyclic, so $\sur{Q}/\sur{M}\sur{N}$ is cyclic. But $\sur{M}\sur{N}$ is a central subgroup of~$\sur{Q}$. It follows that $\sur{Q}$ is abelian, so $H\cong\sur{Q}/\sur{M}$ is abelian. Hence $\sur{Q}/\sur{M}$ is non-cyclic, and it has a subgroup $\sur{M}\sur{N}/\sur{M}$ of order $p$ such that the corresponding quotient $\sur{Q}/\sur{M}\sur{N}$ is cyclic. It follows that $\sur{Q}/\sur{M}\cong H\cong C_p\times C_{p^{h-1}}$, for some $h\geq 2$.\spn
\noindent{\bf Step 5:} Assume that $H$ is neither cyclic nor isomorphic to $C_p\times C_{p^{h-1}}$, for some $h\geq 2$. Then the matrix of the bilinear form $\ncalv{\;}{\,}_{\F,Q}$ is congruent to the identity matrix modulo $p$. In particular, it is non-singular, and the $\F\Aut(Q)$-module $\mathcal{V}_H(Q)=\mathcal{Q}_H^\sharp(Q)$ is isomorphic to the permutation module on the set $\mathcal{N}_H(Q)$.\par
It follows that the $\Aut(P)$-module $\tilde{e}_P^PS_{H,\F}(P)$ is isomorphic to the permutation module on the set of normal subgroups $M$ of $P$ such that $P/M\cong H$. Going back to Step~1 and to the $p$-elementary subgroup $K=P\times C$ of $G$, it follows that the space $\tilde{e}_K^KS_{H,\F}(K)^{N_G(K)}$ has a basis in one to one correspondence with the $N_G(K)$-orbits of normal subgroups $M$ of $K$ such that $K/M\cong H$. Now the isomorphism~(\ref{direct sum decomposition}) shows that $S_{H,\F}(G)$ has a basis in one to one correspondence with the $G$-conjugacy classes of sections $(K,M)$ of $G$ such that $K$ is $p$-elementary and $K/M\cong H$. This proves the theorem, in the case where $H$ is neither cyclic nor isomorphic to $C_p\times C_{p^{h-1}}$, for some $h\geq 2$. \spn
{\bf Step 6: } Suppose now that $H$ is cyclic, of order $p^h>1$. Assume first that $Q$ is cyclic. Then since $Q\cong E\times H$ for some elementary abelian $p$-group $E$, it follows that $E=\un$, i.e. $Q\cong H$. In this case $\mathcal{N}_H(Q)=\{\un\}$, and $\ncalv{\un}{\un}_{H,\F}=1$. Hence $\mathcal{V}_H(Q)$ is isomorphic to the trivial $\F\Aut(Q)$-module in this case.\par
If $Q$ is non-cyclic, let $M,N\in\mathcal{N}_H(Q)$. Recall that
$$\ncalv{M}{N}_{\F,Q}=m_{Q,M\cap N}\,\fracb{\alpha_{M,N}}{|M:M\cap N|}\mvirg$$
where $\alpha_{M,N}$ is defined in (\ref{alphaMN}).\par
The groups $\sur{M}=M/(M\cap N)$ and $\sur{N}=N/(M\cap N)$ are non-trivial elementary abelian central subgroups of $\sur{Q}=Q/(M\cap N)$, and have a common complement in $\sur{Q}$. As $\sur{M}$ is isomorphic to the subgroup $MN/N$ of the cyclic group $Q/N\cong H$, it follows that $\sur{M}\cong C_p$. Moreover $\sur{M}$ has a complement in~$\sur{Q}$, so $\sur{Q}\cong C_p\times C_{p^h}$. Hence if $Q/\Phi(Q)$ has rank $n$, then $\Phi(Q)(M\cap N)/\Phi(Q)$ has rank $n-2$ since
$$Q/\big(\Phi(Q)(M\cap N)\big)\cong \sur{Q}/\Phi(\sur{Q})\cong C_p\times C_p\mpoint$$
By Equations~\ref{mEF} and~\ref{mPQ}, it follows that
$$m_{Q,M\cap N}=(1-p^{n-2})(1-p^{n-3})\cdots(1-p)\mpoint$$
Moreover since $m=1$ and $\gamma=1$, Equation~\ref{scalar} gives
$$\ncalv{M}{N}_{\F,Q}=-m_{Q,M\cap N}(p^s\!-\!1)(p^{s-1}\!-\!1)\cdots(p\!-\!1)p^{\frac{1}{2}{(1-s)(2+s)}-1}\mpoint$$
Since $0\leq s\leq m=1$, there are two cases: \begin{itemize}
\item If $s=1$, then $\sur{M}$ maps into $\sur{Q}/\sur{N}\Phi(\sur{Q})\cong H/\Phi(H)\cong C_p$, hence $MN=Q$ as above, and $Q/N\cong C_{p^h}\cong \sur{M}\cong C_p$, so $h=1$. In this case
\begin{eqnarray*}
\ncalv{M}{N}_{\F,Q}&=&-(1-p^{n-2})(1-p^{n-3})\cdots(1-p)(p-1)/p\\
&=&(1-p^{n-2})(1-p^{n-3})\cdots(1-p^2)(1-p)^2/p\mpoint
\end{eqnarray*}
\item If $s=0$. Then $\sur{M}\leq \sur{N}\Phi(\sur{Q})$, so $\sur{M}\Phi(\sur{Q})=\sur{N}\Phi(\sur{Q})$. If $h=1$, then $\sur{Q}/\sur{M}\cong C_p$, so $\sur{M}\geq\Phi(\sur{Q})$, and it follows that $\sur{M}=\sur{N}$, a contradiction. Thus $h>1$ in this case. Moreover
$$ \ncalv{M}{N}_{\F,Q}=-(1-p^{n-2})(1-p^{n-3})\cdots(1-p)\mpoint$$
\end{itemize}
So in any case, there is a non-zero rational number $\rho$, depending only on $Q$ (and $H$), such that $\ncalv{M}{N}_{\F,Q}=\rho$ when $\ncalv{M}{N}_{\F,Q}\neq 0$. Moreover $\ncalv{M}{N}_{\F,Q}\neq 0$ if and only if $M\neq N$.\par
So the matrix of the form $\ncalv{\;}{\,}_{\F,P}$ is equal to $\rho J$, where $J$ is a matrix of size $|\mathcal{N}_H(Q)|$, with zero diagonal, and non-diagonal coefficients equal to 1. Hence this matrix is non-singular if and only if $|\mathcal{N}_H(Q)|>1$. \par
But $Q= E\times L$, where $L\cong H$ and $E$ is a non-trivial elementary abelian $p$-group. The elements of $\mathcal{N}_H(Q)$ are exactly the groups 
$$E_\varphi=\{\big(e,\varphi(e)\big)\mid e\in E\}\mvirg$$
where $\varphi$ is a group homomorphism from $E$ to $L$. There are $|E|$ such homomorphisms, hence $|\mathcal{N}_H(Q)|=|E|>1$.\par
It follows that the matrix of the form $\ncalv{\;}{\,}_{\F,Q}$ is non-singular, hence the form $\ncalv{\;}{\,}_{\F,Q}$ is non-degenerate. \par
So either when $Q$ is cyclic, or when it is not, the form $\ncalv{\;}{\,}_{\F,Q}$ is non-degenerate. By the same argument as at the end of Step 4, this proves that $S_{H,\F}(G)$ has a basis in one to one correspondence with the $G$-conjugacy classes of sections $(K,M)$ of $G$ for which $K$ is $p$-elementary and $K/M\cong H$. This proves the theorem in the case where $H$ is cyclic.\spn
{\bf Step 7:} Suppose now that $H\cong C_p\times C_{p^{h-1}}$, for some $h\geq 2$. Note that if $h=2$, then $H$ is elementary abelian, so $Q=P/R\cong E\times H$ is elementary abelian. Since $R\leq \Phi(P)$, this forces $R=\Phi(P)$.\par
Now if $M,N\in \mathcal{N}_H(Q)$, since $\gamma=2$ in this case,
$$\ncalv{M}{N}_{\F,Q}=(-1)^m m_{Q,M\cap N}(p^s\!-\!1)(p^{s-1}\!-\!1)\cdots(p\!-\!1)p^{\frac{1}{2}(m-s)(m+s+1)}\mvirg$$
and moreover
$$m_{Q,M\cap N}=(1-p^{n-2})(1-p^{n-3})\cdots(1-p^{m+1})\mvirg$$
where $n$ is the rank of $Q/\Phi(Q)$,  where $m$ is the rank of the elementary abelian groups $\sur{M}=M/(M\cap N)$ and $\sur{N}=N/(M\cap N)$ of $\sur{Q}=Q/(M\cap N)$, and $s$ is the rank of $M/\big(M\cap N\Phi(Q)\big)\cong \sur{M}/\big(\sur{M}\cap \sur{N}\Phi(\sur{Q})\big)$.
Since the exponent $\frac{1}{2}(m-s)(m+s+1)$ of $p$ is non-negative, it follows that $\ncalv{M}{N}_{\F,Q}$ is an integer. Moreover, if $m>s$, this integer is a multiple of $p$. On the other hand $m=s$ if and only if $\sur{M}\cap \sur{N}\Phi(\sur{Q})=\un$, or equivalently if $M\cap N\Phi(Q)=M\cap N$. Since $M\cap\Phi(Q)=N\cap\Phi(Q)=\un$, this is equivalent to $MN\cap \Phi(Q)=\un$. In this case
\begin{equation}\label{rank 2}\ncalv{M}{N}_{\F,Q}=(1-p^{n-2})(1-p^{n-3})\cdots(1-p)
\end{equation}
is congruent to 1 modulo $p$. It follows that the matrix of the form $\ncalv{\;}{\,}_{\F,Q}$ is congruent modulo $p$ to the incidence matrix of the relation $\sim$ on $\mathcal{N}_H(Q)$ defined by $M\sim N$ if and only if $MN\cap\Phi(Q)=\un$. There are now two cases:\spn
$\bullet$ {\bf Case 1:} Assume first that $h\geq 3$, i.e. that $H$ is not elementary abelian of rank 2. 
\begin{mth}{Lemma} \label{matrice}Let $H=C_p\times C_{p^{h-1}}$, for $h\geq 3$, and $Q=E\times H$, where $E$ is an elementary abelian $p$-group of rank $e$. Let $S$ denote the incidence matrix of the relation $\sim$ on $\mathcal{N}_H(Q)$ defined by
$$M\sim N\Leftrightarrow MN\cap\Phi(Q)=\un\mpoint$$
Then:
\begin{enumerate}
\item if $e=0$, the matrix $S$ is the matrix $(1)$.
\item if $e\geq 1$, the eigenvalues of $S$ are $p^{e+1}-p+1$, $p^e-p+1$, and $1-p$, with respective multiplicities $1$, $p^{e+1}-p$, and $p^{2e}-p^{e+1}+p-1$.
\end{enumerate}
In both cases $M$ is invertible modulo $p$.
\end{mth}
\pf If $e=0$, then $E=\un$ and $Q\cong H$, so $\mathcal{N}_H(Q)$ consists of the trivial subgroup $E$ of $Q$. Since $E\sim E$, Assertion~1 follows.\par
If $e\geq 1$, then $\mathcal{N}_H(Q)$ consists of the subgroups
$$E_{\varphi}=\{\big(x,\varphi(x)\big)\mid x\in E\}\mvirg$$
where $\varphi:E\to H$ is a group homomorphism. Since $E$ is elementary abelian, the image of $\varphi$ is contained in the subgroup $C_p\times C_p$ of $H=C_p\times C_{p^{h-1}}$. So there are group homomorphisms $a,b:E\to C_p$ such that $\varphi=(a,b)$, i.e. $\varphi(x)=\big(a(x),b(x)\big)$ for any $x\in E$.\par
Let $\varphi=(a,b)$ and $\varphi'=(a',b')$ be two group homomorphisms from $E$ to~$H$. Then, with an additive notation
$$E_\varphi E_{\varphi'}=\{\big(x-x',a(x)-a(x'),b(x)-b(x')\big)\mid x,x'\in E\}\leq E\times C_p\times C_p\mpoint$$
The element $\big(x-x',a(x)-a'(x'),b(x)-b'(x')\big)$ is in $\Phi(Q)=\un\times\un\times C_{p^{h-2}}$ if and only if $x=x'$ and $a(x)=a'(x')$. Thus
$$E_\varphi\sim E_{\varphi'}\Leftrightarrow \Ker(a-a')\leq \Ker(b-b')\mpoint$$
Identifying $E$ with the vector space $(\F_p)^e$, and $C_p$ with $\F_p$, the homomorphisms $a,b,a',b'$ become elements of the dual vector space $E^*$, and the condition $\Ker(a-a')\leq \Ker(b-b')$ means that there is a scalar $\lambda\in\F_p$ such that $b-b'=\lambda(a-a')$. Hence the incidence matrix $S$ is the matrix indexed by pairs $\big((a,b),(a',b')\big)$ of pairs of elements of $E^*$, defined by
$$S\big((a,b),(a',b')\big)=\left\{\begin{array}{cl}1 &\hbox{if}\;\;\exists \lambda\in\F_p,\;b-b'=\lambda(a-a')\\0&\hbox{otherwise}\end{array}\right.\mpoint$$ 
Let $T$ be the rectangular matrix indexed by the set of pairs $\big((a,b),(c,\lambda)\big)$, where $a,b,c\in E^*$, and $\lambda\in\F_p$, defined by
$$T\big((a,b),(c,\lambda)\big)=\left\{\begin{array}{cl}1 &\hbox{if}\;\;c=b-\lambda a\\0&\hbox{otherwise}\end{array}\right.\mpoint$$ 
Then for $a,b,a',b'\in E^*$, consider the sum
$$s=\sumb{c\in E^*}{\lambda\in\F_p}T\big((a,b),(c,\lambda)\big)T\big((a',b'),(c,\lambda)\big)\mpoint$$
The non-zero terms in this summation correspond to pairs $(c,\lambda)$ such that $c=b-\lambda a=b'-\lambda a'$. Hence $s$ is equal to the number of $\lambda\in\F_p$ such that $b-\lambda a=b'-\lambda a'$. This is equal to 1 if $a'\neq a$ and if $b'-b$ is a scalar multiple of $a'-a$, to $p$ if $a=a'$ and $b=b'$, and to 0 if $a=a'$ and $b\neq b'$. In other words
$$T\cdot{^tT}=S+(p-1)\Id\mpoint$$
Since $T\cdot{^tT}$ is symmetric, it is diagonalizable over $\C$, with real eigenvalues. Let $\mu$ be an eigenvalue of $T\cdot{^tT}$, and $u$ be a corresponding eigenvector. Then $T\cdot{^tT}u=\mu u$, thus ${^tT}\cdot T\cdot{^tT}u=\mu{^tT}u$. So either ${^tT}u=0$, and then $\mu=0$. And if $\mu\neq 0$, then ${^tT}u$ is an eigenvector of ${^tT}\cdot T$ for the eigenvalue $\mu$. Moreover, the map $u\mapsto {^tT}u$ is an injection of the $\mu$-eigenspace of $T\cdot{^tT}$ into the $\mu$-eigenspace of ${^tT}\cdot T$. The same argument applied to ${^tT}\cdot T$ instead of $T\cdot{^tT}$ shows that these two matrices have the same non-zero eigenvalues, and the same multiplicities. \par
Now for $c,c'\in E^*$ and $\lambda,\lambda'\in \F_p$
$$
{^tT}\cdot T\big((c,\lambda),(c',\lambda')\big)=\sum_{a,b\in E^*}T\big((a,b),(c,\lambda)\big)T\big((a,b),(c',\lambda')\big)\mpoint$$
The right hand side is the number of pairs $(a,b)$ of elements of $E^*$ such that $c=b-\lambda a$ and $c'=b-\lambda' a$, i.e. the number of elements $a\in E^*$ such that $c+\lambda a=c'+\lambda' a$, or $c-c'=(\lambda-\lambda')a$. This is equal to $1$ if $\lambda\neq\lambda'$, to $|E|$ if $\lambda=\lambda'$ and $c=c'$, and to 0 if $\lambda=\lambda'$ and $c\neq c'$. Hence the matrix ${^tT}\cdot T$
is a block matrix of the following form
$${^tT}\cdot T=\left(\begin{array}{cccc}|E|\Id&\Omega&\cdots&\Omega\\
\Omega&|E|\Id&\cdots&\Omega\\
\vdots&\vdots&\ddots&\vdots\\
\Omega&\Omega&\ldots&|E|\Id
\end{array}\right)\mvirg$$
where all the $p^2$-blocks are square matrices of size $|E|$, and $\Omega$ is a matrix with all entries equal to 1. Let $\mu$ be an eigenvalue of this matrix, and 
$$v=\left(\begin{array}{c}X_1\\X_2\\\vdots\\X_p\end{array}\right)$$ 
be a corresponding eigenvector, where $X_1,\ldots,X_p$ are column vectors of size~$|E|$. Equivalently, for each $i\in\{1,\ldots,p\}$
$$|E|X_i+\sum_{j\neq i}\Omega X_j=\mu X_i\mpoint$$
But $\Omega X=s(X)\omega$ for any column vector $X$ of size $|E|$, where $s(X)$ denotes the sum of the entries of $X$, and $\omega$ is a column vector of size $|E|$ with all entries equal to 1. Setting $\sigma=\sum_{j=1}^ps(X_j)$, this gives, since $|E|=p^e$
$$p^eX_i+\big(\sigma-s(X_i)\big)\omega=\mu X_i\mpoint$$
Hence if $\mu\neq p^e$, the vector $X_i$ is a multiple of $\omega$, i.e. $X_i=\alpha_i\omega$ for some scalar $\alpha_i$. Then $s(X_i)=\alpha_i p^e$, thus $\sigma=\tau p^e$, where $\tau=\sum_{j=1}^p\alpha_j$. Finally
$$p^e\alpha_i+\big(\tau-\alpha_i)p^e=\tau p^e=\mu\alpha_i\mpoint$$
Thus if $\mu\neq 0$, all the $\alpha_i$'s are equal to $\alpha$, say, and then $\tau=p\alpha$, thus $\mu=p^{e+1}$. Conversely, if $X_i=\omega$ for all $i$, then $v$ is an eigenvector of ${^tT}\cdot T$ with eigenvalue $p^{e+1}$. So $p^{e+1}$ is an eigenvalue of ${^tT}\cdot T$, with multiplicity 1.\par
If $\mu=0$, then the vector $v$ corresponding to $X_i=\alpha_i\omega$ for $i\in\{1,\ldots p\}$ is in the kernel of ${^tT}\cdot T$ if and only if $\sum_{j=1}^p\alpha_i=0$. Hence $0$ is an eigenvalue of ${^tT}\cdot T$ with multiplicity $p-1$. \par
Finally, if $\mu=p^e$, then $s(X_i)=\sigma$ for $i\in\{1,\ldots p\}$, hence $\sigma=p\sigma=0$. The vector $v$ is in the $p^e$-eigenspace of ${^tT}\cdot T$ if and only if $s(X_i)=0$ for all~$i$. Thus $p^e$ is an eigenvalue of ${^tT}\cdot T$, with multiplicity $p(p^e-1)$.\par
It follows that $T\cdot {^tT}$ has eigenvalues $p^{e+1}$, $p^e$, and 0, with respective multiplicities $1$, $p^{e+1}-p$, and $p^{2e}-p^{e+1}+p-1$. This completes the proof, since $S=T\cdot{^tT}-(p-1)\Id$.\findemo
Lemma~\ref{matrice} shows that the form $\ncalv{\;}{\,}_{\F,Q}$ is non-degenerate whenever $H$ is a quotient of $Q$. By the argument of the end of Step 4, or the end of Step 6, this shows that $S_{H,\F}(G)$ has a basis in bijection with the $G$-conjugacy classes of sections $(T,S)$ of $G$ such that $T$ is $p$-elementary and $T/S\cong H$.\spn
$\bullet$ {\bf Case 2:} Suppose finally that $H=C_p\times C_p$. As observed earlier, in this case, if $R$ is a normal subgroup of $P$ contained in $\Phi(P)$ such that $P/R\cong E\times H$ for some elementary abelian $p$-group $E$, then in fact $R=\Phi(P)$. The group $Q=P/R$ is elementary abelian, and decomposes as $Q=E\times L$, where $L\cong H$ is elementary abelian of rank~2. The set $\mathcal{N}_H(Q)$ is the set of complements $M$ of $L$ in $Q$, and $MN\cap\Phi(Q)=\un$ for any $M,N\in\mathcal{N}_H(Q)$. Equation~\ref{rank 2} shows that 
$$\ncalv{M}{N}_{\F,Q}=(1-p^{n-2})(1-p^{n-3})\cdots(1-p)\mvirg$$
where $n$ is the rank of $P/\Phi(P)$. {\em This is non-zero, and does not depend on $M,N\in\mathcal{N}_H(Q)$}.  Hence the form $\ncalv{\;}{\,}_{\F,Q}$ has rank 1 in this case. Thus $\tilde{e}_P^PS_{H,\F}(P)$ is one dimensional if $P$ is non-cyclic, and it is zero otherwise. Saying that $P$ is non-cyclic is equivalent to saying that the $p$-elementary group $K=P\times C$ of Step~1 is non-cyclic. Hence $S_{H,\F}(G)$ has a basis in bijection with the conjugacy classes of non-cyclic $p$-elementary subgroups of~$G$. This completes the proof of Theorem~\ref{the theorem}.\findemo
\begin{rem}{Remark} As $C_p\times C_p$ is a $B$-group, Case 2 above also follows from Proposition 11 of~\cite{doublact}: for a $B$-group $H$ and a finite group $G$, the dimension of $S_{H,\F}(G)$ is equal to the number of conjugacy classes of subgroups $K$ of~$G$ such that $\beta(K)\cong H$. Now by Lemma~\ref{p-elementary}, if $\beta(K)\cong C_p\times C_p$, then $K$ is $p$-elementary, and non cyclic (for otherwise $\beta(K)=\un$). Conversely, if $K$ is $p$-elementary and non cyclic, then $\beta(K)$ is a non trivial $p$-group, and also a $B$-group, hence $\beta(K)\cong C_p\times C_p$.
\end{rem}
\section{A Green biset functor for $p$-elementary groups}
The following theorem is closely related to Theorem~\ref{the theorem}. In particular, it yields an alternative proof of its Assertions 1 and 2. We refer to Section 8.5 of~\cite{bisetfunctorsMSC} for the basic definitions on Green biset functors. 
\begin{mth}{Theorem} \label{un foncteur}Let $p$ be a prime number. 
\begin{enumerate}
\item For a finite group $G$, let $\mathcal{E}l_p(G)$ denote the set of $p$-elementary subgroups of $G$. Set 
$$F_p(G)=\{u\in B(G)\mid\forall H\in \mathcal{E}l_p(G),\;\Res_H^Gu=0\}\mpoint$$
Then the assignment $G\mapsto F_p(G)$ is a biset subfuntor of the Burnside functor $B$, and the quotient functor
$$E_p =B/F_p$$
is a Green biset functor (over $\Z)$.
\item For a finite group $G$, the evaluation $E_p(G)$ is a free abelian group of rank equal to the number of conjugacy classes of $p$-elementary subgroups of $G$.
\item Let $\F$ be a field of characteristic 0. Then the biset functor $\F E_p=\F\otimes_\Z E_p$ has a unique non zero proper subfunctor $I$, isomorphic to $S_{(C_p)^2,\F}$, and the quotient $\F E_p/I$ is isomorphic to $S_{\un,\F}\cong \F R_\Q$. In other words there is a non split short exact sequence
\begin{equation}\label{exact}
0\to S_{(C_p)^2,\F}\to \F E_p\to S_{\un,\F}\to 0
\end{equation}
of biset functors over $\F$.
\end{enumerate}
\end{mth}
\pf Let $\F B=\F\otimes _\Z B$ be the Burnside functor over $\F$. If we forget the $\F$-structure on $\F B$, we get an inclusion $B\to\F B$ of biset functors over~$\Z$. In particular, for each finite group $G$, we get an inclusion 
$$f_p:F_p(G)\to \F B(G)\mpoint$$
Now saying that $u\in B(G)$ lies in $F_p(G)$ amounts to saying that the restriction of $f_p(u)$ to any $p$-elementary subgroup of $G$ is equal to 0. Since any subgroup of a $p$-elementary group is again $p$-elementary, this amounts to saying that $|f_p(u)^H|=0$ for any $H\in\mathcal{E}l_p(G)$. In other words $f_p(u)$ is a linear combination of idempotents $e_K^G$ of $\F B(G)$, where $K$ is a subgroup of~$G$ which is {\em not} $p$-elementary. By Lemma~\ref{p-elementary}, we get that $u\in\F_p(G)$ if and only if $f_p(u)$ is a linear combination of idempotents $e_K^G$, for subgroups $K$ such that $\beta(K)$ is not a $p$-group, that is $\beta(K)$ is non trivial and not isomorphic to $(C_p)^2$.\par
Let $\mathcal{G}_p$ be the class of $B$-groups which are non trivial, and not isomorphic to $(C_p)^2$. Then $\mathcal{G}_p$ is a {\em closed} class of $B$-groups (\cite{bisetfunctorsMSC}, Definition 5.4.13), that is, if a $B$-group $L$ admits a quotient in $\mathcal{G}_p$, then actually $L\in \mathcal{G}_p$ (this is because the only quotient $B$-groups of $(C_p)^2$ are the trivial group and $(C_p)^2$, up to isomorphism). By Theorem~5.4.14 of~\cite{bisetfunctorsMSC}, this closed class $\mathcal{G}_p$ is associated to a subfunctor $N_p$ of the Burnside functor $\F B$, defined for a finite group $G$ by
$$N_p(G)=\sum_{\substack{K\leq G\\\beta(K)\in \mathcal{G}_p}}\F e_K^G\mpoint$$
This shows that $f_p\big(F_p(G)\big)=f_p\big(B(G)\big)\cap N_p(G)$, and since $N_p$ is a biset subfunctor of $\F B$, it follows that $F_p$ is a biset subfunctor of $B$. As biset subfunctors of $B$ are also ideals of the Green biset functor $B$ (see Lemma~2.5.8, Assertion 4 in~\cite{bisetfunctorsMSC}), we get that $F_p$ is an ideal of $B$. It follows that the quotient $E_p=B/F_p$ is a Green biset functor. This completes the proof of Assertion~1.\par
Moreover, the ghost map
$$\Phi_G:B(G)\to\prod_{\substack{K\leq G\\{\rm mod.}\, G}}\Z$$
sending $u\in B(G)$ to the sequence $|u^K|$, is injective by Burnside's theorem.
 The above discussion shows that $\Phi$ induces an injective map
$$E_p(G)=B(G)/F_p(G)\to \prod_{\substack{K\in \mathcal{E}l_p(G)\\{\rm mod.}\, G}}\Z\mvirg$$
which becomes an isomorphism after tensoring with $\F$. Assertion 2 follows.\par
Finally, it follows from Theorem~5.4.14 of~\cite{bisetfunctorsMSC} that the lattice $[0,\F E_p]$ of biset subfunctors of $\F E_p$ is isomorphic to the set of closed classes of $B$-groups which contain $\mathcal{G}_p$. There are exactly three such classes: the class $\mathcal{G}_p$, the class of non-trivial $B$-groups, and the class of all $B$-groups. So $[0,\F E_p]$ is a totally ordered set of cardinality 3. Hence $\F E_p$ admits a unique non zero proper subfunctor $I$. The quotient $\F E_p/I$ is the unique simple quotient of $\F B$, hence it is isomorphic to $S_{\un,\F}\cong\F R_\Q$. Now $I$ is a simple biset functor, which is a subquotient of $\F B$. By Proposition~5.5.1 of~\cite{bisetfunctorsMSC}, it follows that $I\cong S_{H,\F}$ for some $B$-group $H$. Since the group $K=(C_p)^2$ has a unique non cyclic subgroup, it follows that $I(K)$ is one dimensional, and a trivial $\F\Out(K)$-module. Moreover $K$ is a group of minimal order such that $I(K)\neq\zero$. Hence $H\cong K$, and $I\cong S_{(C_p)^2,\F}$. This completes the proof of Assertion 3, and the proof of Theorem~\ref{un foncteur}.\findemo
\begin{rem}{Remark} One can show that the exact sequence~(\ref{exact}) is essentially unique as a non split exact sequence in the category $\mathcal{F}_\F$ of biset functors over~$\F$: more precisely, one can show that $\Ext_{\mathcal{F}_\F}^1(S_{\un,\F},S_{(C_p)^2,\F})\cong\F$.
\end{rem}
\begin{rem}{Remark} If $G$ is a $p$-group (or even if $G$ is $p$-elementary), then $F_p(G)=\zero$, so $E_p(G)\cong B(G)$. So if we restrict the exact sequence~(\ref{exact}) to finite $p$-groups, we get an exact sequence
$$ 0\to S_{(C_p)^2,\F}\to \F B\to \F R_\Q\to 0$$
of $p$-biset functors over $\F$. This (restricted) exact sequence was introduced in~\cite{both}, where is was shown that for a finite $p$-group $P$, the evaluation $S_{(C_p)^2,\F}(P)$ is isomorphic to $\F D(P)$, where $D(P)$ is the Dade group of endopermutation modules. It was also shown that the dimension of $S_{(C_p)^2,\F}(P)$ is equal to the number of conjugacy classes of non cyclic subgroups of $P$, which is also the number of conjugacy classes of non-cyclic $p$-elementary subgroups of $P$. So this agrees with Assertion 2 of Theorem~\ref{the theorem}.
\end{rem}
\begin{rem}{Remark} For a finite group $G$, let $M_p(G)$ be the $\Z$-submodule of $B(G)$ generated by the classes of the transitive $G$-sets $G/H$, where $H$ is a $p$-elementary subgroup of $G$. One can check easily that $M_p(G)\cap F_p(G)=\zero$, so comparing ranks, one might hope that $B(G)=M_p(G)\oplus F_p(G)$. This is false in general: for $p=2$, when $G$ is the symmetric group $S_3$, there are three $p$-elementary subgroups in $G$, up to conjugation, namely the proper subgroups of $G$ (that is the trivial group, the alternating subgroup $A=A_3$, and the subgroup $C$ of order 2). Hence if $B(G)=M_p(G)\oplus F_p(G)$, then in particular $G/G\in M_p(G)\oplus F_p(G)$ so there exist integers $a,b,c$ such that the element $u=G/G-(b\,G/\un+a\,G/A+c\,G/C)$ is in $F_p(G)$. Taking fixed points by $A$ then gives $|u^A|=0=1-2a$, a contradiction. One can show more precisely that $M_p(G)\oplus F_p(G)$ has index 2 in $B(G)$ in this case.
\end{rem}


\vspace{3ex}
\begin{flushleft}
Serge Bouc\\
LAMFA-CNRS UMR 7352\\
Universit\'e de Picardie-Jules Verne\\
33, rue St Leu, 80039 Amiens Cedex 01\\
France\\
{\small\tt serge.bouc@u-picardie.fr}\\
{\small\tt http://www.lamfa.u-picardie.fr/bouc/}
\end{flushleft}

\end{document}